\documentclass[11pt]{article}
\usepackage{epic,latexsym}
\usepackage{amssymb,amsmath,amsthm}
\usepackage{color}
\usepackage{tikz}
\usepackage{tabu}
\usepackage{multirow}

\newtheorem{thm}{Theorem}
\newtheorem{lem}[thm]{Lemma}

\newtheorem{cor}{Corollary}

\newenvironment{unnumbered}[1]{\trivlist \item [\hskip \labelsep {\bf
#1}]\ignorespaces\it}{\endtrivlist}

\newcommand{\smallqed}{{\tiny ($\Box$)}}

\def \nmr {\begin{enumerate}}
\def \enmr {\end{enumerate}}

\def \tmz {\begin{itemize}}
\def \etmz {\end{itemize}}

\newcommand{\Proof}{\noindent\textbf{Proof. }}

\newcommand{\w}{{\rm w}}

\begin{document}

\title{Bounds on the $2$-domination number}
\author{$^{1,2}$Csilla Bujt\'{a}s \enskip and \enskip
$^1$Szil\'ard Jask\'o
\\ \\
$^1$Faculty of Information Technology \\
University of Pannonia, Veszpr\'em, Hungary\\
\small \tt Email: bujtas@dcs.uni-pannon.hu, \enskip
jasko.szilard@uni-pen.hu
   \\
\\
$^2$Alfr\'ed R\'enyi Institute of Mathematics \\
       Hungarian Academy of Sciences, Budapest, Hungary \\
}
\date{}
\maketitle

\begin{abstract}
In a graph $G$, a set $D\subseteq V(G)$ is called $2$-dominating set
if each vertex not in $D$ has at least two neighbors  in $D$.
 The  $2$-domination number $\gamma_2(G)$ is  the minimum cardinality of
 such a set $D$.
We give a method for the construction of 2-dominating sets, which
also yields upper bounds on the 2-domination number in terms of the
 number of vertices, if the minimum degree $\delta(G)$ is fixed. These improve the best earlier bounds
 for any  $6 \le \delta(G) \le 21$. In particular, we prove that $\gamma_2(G)$ is strictly smaller than $n/2$, if $\delta(G) \ge 6$.
 Our proof technique uses
a  weight-assignment to the vertices  where the weights are
 changed during the procedure.
\end{abstract}


\section{Introduction}

We study the  graph invariant $\gamma_2(G)$, called 2-domination
number, which is in close connection with the fault-tolerance of
networks. Our main contributions are upper bounds on $\gamma_2(G)$
in terms of the number of vertices, when the minimum degree
$\delta(G)$ is fixed. The earlier upper bounds of this type are
tight for $\delta(G) \le 4$, here we establish improvements for the
range of $6 \le \delta(G) \le 21$. Our approach is based on a
weight-assignment to the vertices, where the weights are changed
according to some rules during a 2-domination procedure.

\subsection{Basic terminology}

Given a simple undirected graph $G$, we denote by $V(G)$ and $E(G)$
the set of its vertices and edges, respectively. The \emph{open
neighborhood} of a vertex $v \in V(G)$ is defined as $N(v)=\{u \in
V(G)\mid uv \in E(G)\}$, while the \emph{closed neighborhood} of $v$
is $N[v]=N(v)\cup \{v\}$. Then, the {\it degree} $d(v)$ is equal to
$|N(v)|$ and the \emph{minimum degree} of  $G$ is  the smallest
vertex degree $\delta(G)=\min\{d(v)\mid v \in V(G)\}$. We say that a
vertex $v$ \emph{dominates} itself and its neighbors, that is
exactly the vertices  contained in $N[v]$. A set $D\subseteq V(G)$
is a \emph{dominating set} if each vertex of $G$ is dominated or
equivalently, if the closed neighborhood of $D$, defined as $N[D]=
\bigcup_{v \in D} N[v]$, equals $V(G)$. The \emph{domination number}
$\gamma(G)$  is the minimum cardinality of such a set $D$.
Domination theory has a rich literature, for results and references
see the monograph \cite{Fu}.

There are two different natural ways to generalize the notion of
(1-)domination to multiple domination. As defined in \cite{FJ}, a
 \emph{$k$-dominating set} is a  set $D\subseteq V(G)$ such that every vertex
not in $D$ has at least $k$ neighbors in $D$. Moreover, $D$ is a
\emph{$k$-tuple dominating set} if the same condition $|N[v] \cap
D|\ge k$ holds not only for all $v\in V(G)\setminus D$ but for all
$v\in V(G)$. The minimum cardinalities of such sets are the
\emph{$k$-domination number} $\gamma_k(G)$ and the \emph{$k$-tuple
domination number} of $G$, respectively.

\subsection{2-domination and applications}

A sensor network can be modeled as a graph such that the vertices
represent the sensors and two vertices are adjacent if and only if
the corresponding devices can communicate with each other. Then, a
dominating set $D$ of this graph $G$ can be interpreted as a
collection of cluster-heads, as each sensor which does not belong to
$D$ has at least one head within communication distance.

A $k$-dominating set $D$ may represent a dominating set which is
$(k-1)$-fault tolerant. That is, in case of the failure of at most
$(k-1)$ cluster-heads, each remaining vertex is either a head or
keeps in connection with at least one head. The price of this
$k$-fault tolerance might be very high. In the extremal case, when
$k$ is greater than the maximum degree in the network, the only
$k$-dominating set is the entire vertex set. But for the usual cases
arising in practice, 2-domination  might be enough and it does not
require extremely many heads.

Note that $k$-tuple domination might need much more vertices
(cluster-heads) than $k$-domination. As proved in \cite{appl}, for
each real number $\alpha >1$  and  each natural number $n$ large
enough, there exists a graph $G$ on $n$ vertices such that its
$k$-tuple domination number is at least $\frac{k}{\alpha}$ times
larger than its $k$-domination number. There surely exist some
practical problems where $k$-tuple domination is needed, but for
many problems arising $k$-domination seems to be sufficient. Indeed,
if a cluster-head  fails and is deleted from the network, we may not
need further heads to supervise it. This motivates our work on the
2-domination number $\gamma_2$.

Another potential application of our results in sensor networks
concerns the data collection problem. Here, each sensor has two
capabilities: either measures and reports, or receives and collects
data. Only one position from those two can be active at the same
time. After deploying, the organization process determines exactly
which sensors supply the measuring and the collector function in the
given network. Since it is a natural condition that every
measurement should be saved in at least two different devices, the
set of collector sensors should form a 2-dominating set in the
network.

We  mention shortly that many further kinds of application exist.
For example a facility location problem may require that each region
is either served by its own facility or has at least two neighboring
regions with such a service \cite{HC}. In this context, facility
location may also mean allocation of a camera system, or that of
ambulance service centers.

\subsection{Upper bounds on the 2-domination number}

Although this subject attracts much attention   (see the recent
survey \cite{CFHV} for  results and references)  and it seems very
natural to give upper bounds for $\gamma_2$ in terms of the minimum
degree, there are not too many  results of this type. The following
general upper bounds are known. (As usual, $n$ denotes the order of
the graph, that is the number of its vertices.)
\begin{itemize}
\item If the minimum degree $\delta(G)$ is $0$ or $1$, then $\gamma_2(G)$
can be equal to $n$.
\item If  $\delta(G)= 2$ then $\gamma_2(G) \le
\frac{2}{3}\;n$. This statement follows from a general upper bound
on $\gamma_k(G)$ proved in \cite{FHV}. The bound is tight for graphs
 each component of which is a $K_3$.
\item If $\delta(G)\ge 3$ then $\gamma_2(G) \le \frac{1}{2}\;n$. The
general theorem, from which the bound follows, was established in
\cite{CR}. Note that a 2-dominating set of cardinality at most $n/2$
can be constructed by a simple algorithm. We divide the vertex set
into two parts and then in each step, a vertex  which has more
neighbors in its own part than in the other one, is moved into the
other part. If the minimum degree is at least 3, this procedure
results in two disjoint 2-dominating sets. Note that for
$\delta(G)=3$ and $4$ the bound is tight. For example, it is easy to
check that  $\gamma_2(K_4)=2$ and  $\gamma_2(K_4 \Box K_2)=
4$.\footnote{The Cartesian product $K_4 \Box K_2$ is the graph of
order 8 which consists of two copies of $K_4$ with a matching
between them. Note that  $\gamma_3(K_4 \Box K_2)$ also equals 4.}
\item For every graph $G$ of minimum degree $\delta \ge 0$,
$$\gamma_2(G) \le \frac{2\ln(\delta +1)+1}{\delta +1}\; n.$$
This upper bound was obtained in \cite{HV}  using probabilistic
method and it is a   strong result when $\delta$ is really high. On
the other hand, it gives an upper bound better than $0.5\; n$ only
if $\delta(G) \ge 11$.
\end{itemize}
In this paper we present a method which can be used to improve the
existing upper bounds when the minimum degree  is in the ``middle''
range. Particularly, we show that
  if $\delta(G)\ge 6$ then $\gamma_2(G)$ is strictly smaller than $n/2$;
  $\delta(G)=7$ implies $\gamma_2(G) < 0.467 \; n$;
  $\delta(G)=8$ implies $\gamma_2(G) <  0.441 \; n$;
  and $\gamma_2(G) < 0.418 \; n$ holds for every graph whose minimum
   degree is at least 9.

The paper is organized as follows. In Section~\ref{sect2}, we state
our main theorem and its corollaries which are the new upper bounds
for specified minimum degrees. In Section~\ref{sec-proof} our main
theorem is proved. Finally, we make some remarks on the algorithmic
aspects of our results.

\section{Our results} \label{sect2}

To avoid the repetition of the analogous argumentations for
different minimum degrees, we will state our theorem in a general
form which is quite technical. Then, the upper bounds will follow as
easy consequences. First, we introduce a set of conditions which
will be referred to in our main theorem. We assume that $d\ge 4$
holds.

\begin{align}
s> a \ge y_{d+1} \ge y_d \ge \dots \ge y_0 \ge b_0 =0\\
 0 \le b_{d+1}-b_d \le  b_d-b_{d-1} \le  \dots \le b_2-b_1 \le b_1\\
 0\le y_{d+1}-b_{d+1} \le y_d-b_d  \le \dots \le y_1-b_1 \le y_0\\
  y_{d+1} \le a-\frac{s-a}{d+2}\\
 y_d \le a-\frac{s-a}{d+1}\\
 b_{d+1}\le a-\frac{s-a}{d+3}-\frac{s-a}{d+1}\\
  b_d\le a-\frac{s-a}{d+2}-\frac{s-a}{d+1}\\
  b_{d-1} \le a-2\cdot \frac{s-a}{d+1}\\
 a+d(a-y_{d-1}) &\ge s\\
 a+d(y_{d+1}-b_d)&\ge s\\ 
 y_{d+1}+(d+1)(a-y_{d-2}) &\ge s\\
 y_{d+1}+(d+1)(y_{d+1}-b_d) &\ge s\\ 
 a+(d-1)(a-y_{d-2}) &\ge s\\
 a+(d-1)(y_d-b_{d-1})&\ge s\\ 
 y_d+d(a-y_{d-3}) &\ge s\\
 y_d+d(y_d-b_{d-1}) &\ge s\\ 
 b_{d+1}+(d+1)(a-y_{d-2})&\ge s\\
 b_{d+1}+(d+1)(y_{d-1}-b_{d-1})&\ge s\\   
  a+(d-1)(b_3-b_2)+(y_2-b_1)+(d-3)(b_3-b_2)&\ge s\\
 y_2+(d-3)(b_3-b_2)+2(y_2-b_1+(d-3)(b_3-b_2))&\ge s\\ 
  a+0,5b_3+(d-1.5)(b_3-b_2)+(a-y_1)&\ge s\\
 a+0,5b_3+(d-1.5)(b_3-b_2)+(y_1-b_1+(d-3)(b_3-b_2))&\ge s\\
 1.5a+ 0.5b_3 - 0.5y_0 + (d-2)(b_2-b_1) &\ge s\\
  a+ 0.5b_3 + 0.5y_1 - 0.5b_1 + (d-2)(b_2-b_1)+0.5(d-3)(b_3-b_2) &\ge s\\
 1.5a +b_3+0.5(d-3)(b_3-b_1) + 0.5(d-2)(b_3-b_2) &\ge s \\
 a +b_3+ 0.5y_1-0.5b_1 +0.5(d-3)(b_3-b_1) + 0.5(d-4)(b_3-b_2) &\ge s
 \end{align}
 \begin{align}
   b_3+3(a-y_0)&\ge s\\
 b_3+3(y_1-b_1+(d-3)(b_3-b_2))&\ge s\\ 
 2y_1+2(d-2)(b_2-b_1)&\ge s \\
 a+0,5b_2+(d-1.5)(b_2-b_1)+0,5(a-y_1)&\ge s\\
  a+0,5b_2+(d-1.5)(b_2-b_1)+0,5(y_1-b_1+(d-3)(b_2-b_1))&\ge s\\
  a+0.5(d-1)  b_2&\ge s\\
       b_2+2y_0+2(d-2)(b_2-b_1)&\ge s\\
 b_2+y_0+(d-2)(b_2-b_1)+(a-y_0)&\ge s\\
 y_0+(d-1)b_1&\ge s
 \end{align}

 For every  $2 \le i \le d-2$:  
 \begin{align}
a+ (d-i)(b_{i+2}-b_{i+1}) +i(a-y_{i-1}) &\ge s \\
 a+(d-i)(b_{i+2}-b_{i+1})+ i(y_{i+1}-b_{i}+(d-i-2)(b_{i+2}-b_{i+1}))&\ge s\\
   y_{i+1}+(d-i-2)(b_{i+2}-b_{i+1})+ (i+1)(a-y_{i-2})&\ge s\\
 y_{i+1}+(d-i-2)(b_{i+2}-b_{i+1})+ (i+1)(y_{i+1}-b_{i}+(d-i-2)(b_{i+2}-b_{i+1}))&\ge s\\ 
 b_{i+2}+(i+2)(a-y_{i-1})&\ge s\\
 b_{i+2}+(i+2)(y_{i}-b_{i}+(d-i-2)(b_{i+2}-b_{i+1}))&\ge s 
  \end{align}

\bigskip
Now we are in a position to state our main theorem. Its proof will
be given in Section~\ref{sec-proof}.

\begin{thm} \label{th1}
Assume that $G$ is a graph of order $n$ and with minimum degree
$\delta(G)=d \ge 6$. If
 $a$, $y_0, \dots, y_{d+1}$, $b_0, \dots, b_{d+1}$ are
nonnegative numbers and $s$ is a positive number such that
conditions $(1)$--$(35)$, and also for every $2\le i \le d-2$ the
inequalities $(36)$--$(41)$ are satisfied, then
$$\gamma_2(G) \le \frac{a}{s}\;n.$$
\end{thm}

If we fix an integer $d$, set $s=1$, and want to minimize $a$ under
the conditions given in Theorem~\ref{th1}, we have a linear
programming problem. The solution $a^*$ of this LP-problem gives an
upper bound on $\frac{\gamma_2(G)}{n}$ which holds for every graph
with $\delta(G) \ge d$. In Table~\ref{Atable11}, we summarize these
upper bounds for several values of $d$.

\begin{table}[]
 \label{Atable11}
 \centering
\begin{tabu}{|[2pt]l|[2pt]c|c|c|c|c|c|[2pt]}
        \tabucline[2pt]{-}
    ~$\delta$ &  $6$ & $7$ & $8$ & $9$ & $10$ & $11$ \\ \hline 
    Our result & $0.49754$ &  $0.46682$ & $0.44016$ & $0.41702$ & $0.39679$ & $0.37957$    \\ \hline
    Earlier best bound &  $0.5$  & $0.5$ & $0.5$ & $0.5$ & $0.5$ & $0.49749$  \\  \tabucline[2pt]{-}
      ~$\delta$ &  $12$ & $13$ & $14$ & $15$ & $16$ & $17$  \\ \hline 
    Our result & $0.36459$ &  $0.35117$ & $0.33914$ & $0.33385$ & $0.33052$ & $0.32762$   \\ \hline
    Earlier best bound & $0.47154$ & $0.44844$ & $0.42775$ & $0.40908$ & $0.39215$ & $0.37671$   \\ \tabucline[2pt]{-}
      ~$\delta$ &  $18$ & $19$ & $20$ & $21$ & $22$ & $23$  \\ \hline 
    Our result & $0.32505$ &  $0.32277$ & $0.32074$ & $0.31891$ & $0.31726$ & $0.31574$   \\ \hline
    Earlier best bound & $0.36258$ & $0.34958$ & $0.33758$ & $0.32646$ & $0.31613$ & $0.30651$   \\ \tabucline[2pt]{-}
      ~$\delta$ &  $24$ & $25$ & $26$ & $27$ & $30$ & $40$  \\ \hline 
    Our result & $0.31436$ &  $0.31309$ & $0.31192$ & $0.31084$ & $0.30803$ & $0.30178$   \\ \hline
    Earlier best bound & $0.29752$ & $0.28909$ & $0.28118$ & $0.27373$ & $0.25381$ & $0.20555$   \\ \tabucline[2pt]{-}
      ~$\delta$ &  $50$ & $60$ & $70$ & $80$ & $90$ & $100$  \\ \hline 
    Our result & $0.29806$ &  $0.29560$ & $0.29385$ & $0.29254$ & $0.29152$ & $0.29071$   \\ \hline
    Earlier best bound & $0.17380$ & $0.15118$ & $0.13416$ & $0.12086$ & $0.11013$ & $0.10129$   \\ \tabucline[2pt]{-}

\end{tabu}
\caption{Comparison of  our results and earlier best
upper bounds on $\frac{\gamma_2(G)}{n}$, if the minimum degree
$\delta$ is fixed.}
  \label{Atable11}
\end{table}


The following consequences for $d=6,7,8,9$ can be directly obtained
by using the integer values given for the variables $s,a, y_0,
\dots, y_{d+1},b_0, \dots, b_{d+1}$ in Table~\ref{Atable22}.
Substituting them into the conditions (1)--(41) of
Theorem~\ref{th1}, one can check that all inequalities are
satisfied. This yields the following upper bounds on the
2-domination number.


\begin{table}[h]
 \centering
\begin{tabu}{|[2pt]l|[2pt]c|c|c|c|[2pt]}
        \tabucline[2pt]{-}
        ~        & $\delta=6$ & $\delta=7$ & $\delta=8$ & $\delta=9$ \\ \tabucline[2pt]{-}
         $a$  & $502562162340$ & $9858456650$ & $215321625855$ & $93641183816180$   \\ \hline
         $s$  & $1010109434040$ & $21118330730$ & $489195209055$  & $224551068595700$  \\ \hline
        $y_{10}$ & $-$ & $-$ & $-$  & $78747157548500$   \\ \hline
        $y_9$ & $-$ & $-$ & $180637395519$  & $78747157548500$   \\ \hline
        $y_8$ & $-$     & $8265018290$ & $180637395519$  & $78747157548500$  \\ \hline
        $y_7$ & $422846061750$     & $8265018290$ & $180637395519$  & $77277448218740$   \\ \hline
        $y_6$ & $422846061750$     & $8093880725$ & $176196828255$  & $75612599739380$  \\ \hline
        $y_5$ & $409645123200$     & $7981810970$ & $170236790715$  & $73000318746740$  \\ \hline
        $y_4$ & $401052708000$     & $7754608778$ & $164408232975$  & $69343125357044$  \\ \hline
        $y_3$ & $387969820875$     & $7321226150$ & $153359038875$  & $64634747985500$  \\ \hline
        $y_2$ & $357968691360$     & $6598921770$ & $138571857655$  & $57524154844772$  \\ \hline
        $y_1$ & $296456709780$     & $5196793700$ & $105895928425$  & $43483590947181$  \\ \hline
        $y_0$ & $254021681340$     & $4492799990$ & $87943795415$  & $33987088151324$  \\ \hline
        $b_{10}$ & $-$     & $-$ & $-$   & $64166766443780$ \\ \hline
        $b_9$ & $-$     & $-$ & $146353194015$   & $64166766443780$ \\ \hline
        $b_8$ & $-$     & $6656464850$ & $146353194015$  & $61811868322820$  \\ \hline
        $b_7$ & $338254849800$     & $6656464850$ & $139847385195$  & $59456970201860$  \\ \hline
        $b_6$ & $338254849800$     & $6286147490$ & $133341576375$  & $57102072080900$  \\ \hline
        $b_5$ & $313665896880$     & $5915830130$ & $126835767555$   & $54125243789540$ \\ \hline
        $b_4$ & $289076943960$     & $5545512770$ & $118110911835$   & $50365444145324$ \\ \hline
        $b_3$ & $264487991040$     & $5021360750$ & $107061717735$   & $45588781601132$ \\ \hline
        $b_2$ & $226888474680$     & $4278173340$ & $89997559750$   & $37861061453138$ \\ \hline
        $b_1$ & $151217550540$     & $2770921790$ & $57321630520$   & $23820497555547$ \\  \tabucline[2pt]{-}
\end{tabu}
\caption{Weights assigned to the vertices for graphs of minimum
degree $\delta=6,7,8$ and $9$.}  \label{Atable22}
\end{table}

\begin{cor} Let $G$ be a graph of order $n$.
\begin{itemize}
\item[$(i)$] If $\delta(G) = 6$ then $\gamma_2(G) \le
\frac{456883}{918298} \; n < 0.498 n$.
\item[$(ii)$] If $\delta(G) = 7$ then $\gamma_2(G) \le \frac{140835095}{301690439} \; n < 0.467 n$.
\item[$(iii)$] If $\delta(G) = 8$ then $\gamma_2(G) \le \frac{292954593}{665571713} \; n < 0.441 n$.
\item[$(iv)$] If $\delta(G) \ge 9$ then $\gamma_2(G) \le
\frac{60805963517}{145812382205}\; n < 0.418 n$. 
\end{itemize}
\end{cor}

\section{Proof of Theorem~\ref{th1}} \label{sec-proof}

To prove Theorem~\ref{th1} we apply an algorithmic approach, where
weights are assigned to the vertices and these weights change
according to some rules during the greedy 2-domination procedure.
  A similar proof technique was introduced in
\cite{bu-2013}, later it was used in \cite{bu-2014, bujtas-2014,
sch-2016} for obtaining upper bounds on the game domination number
(see \cite{BKR} for the definition) and in \cite{HKR-2016, HR-2016}
for proving bounds on the game total domination number
\cite{HKR-2015}. Based on this approach we also obtained
improvements for the upper bounds on the domination number
\cite{BK}, and in the conference paper \cite{bu-vocal} we presented
a preliminary version of this algorithm to estimate the 2-domination
number of graphs of minimum degree 8.

\subsection{Selection procedure with changing weights}

 Throughout, we assume that a graph $G$ is given with $\delta (G) \ge d \ge 6$.
 We will consider an algorithm in which the vertices of the 2-dominating set  are
 selected one-by-one. A step in the algorithm means that one vertex is selected (or chosen) and put into the set $D$ which was empty at the beginning of the process.
  Hence, after any step of the procedure, $D$ denotes the set of vertices
  chosen up to this point. We make difference between the following four main types of
  vertices:
  \begin{itemize}
  \item A vertex $v$ is {\it white},   if $v$ is not dominated,  that is
  if $|N[v] \cap D| = 0$.
  \item A vertex $v$ is {\it yellow}, if $|N(v)\cap D|=1$  and $v \notin D$.
  \item A vertex $v$ is {\it blue}, if $|N(v)\cap D|\ge 2$
  and $v \notin D$.
  \item A vertex $v$ is {\it red}, if $v\in D$.
  \end{itemize}
  The sets of  the white, yellow, blue and red vertices are denoted by $W$, $Y$, $B$
  and $R$, respectively.  After any step of the algorithm, we consider
  the graph $G$ together with the set $D$.  Hence,  the current colors of the vertices, that is the partition
  $V(G)=W \cup Y \cup B \cup R$, are also determined. The graph $G$ together with a $D\subseteq V(G)$  will be called {\it colored graph}
  and denoted by $G^D$.
  We define the  {\it WY-degree}
  of a vertex $v$ in $G^D$ to be $\deg_{WY}(v)= |N(v) \cap (W \cup Y)|$.
  The sets $W$, $Y$ and $B$ are partitioned according to the
  WY-degrees of the vertices. For every integer $i \ge 0$ and for $X=W, Y, B$, let
   $X_i=\{v\in X \mid \deg_{WY}(v)=i\}$. Since $R=D$, we may assume
   that red vertices are not selected in any steps of the procedure.

 We distinguish between two types of colored graphs. $G^D$ belongs to
 {\it Type 1} if
 $\max\{i \mid W_{i} \cup Y_{i+1}\neq \emptyset\} \ge d+1$,
 otherwise $G^D$ is of {\it Type 2}.
 Hence, a colored graph is of Type 2 if and only if
 $\deg_{WY}(v) \le d$  for every white vertex $v$ and  $\deg_{WY}(u) \le d+1$ for every yellow
 vertex $u$.

 During the 2-domination algorithm, weights are assigned to the vertices. The weight $\w(v)$ of vertex $v$ is defined with respect to the current
 type of the colored graph
  and to the current color and WY-degree of  $v$.

\begin{table}[h]
 \label{table2}
 \centering
\begin{tabu}{|[2pt]l|[2pt]cl|c|[2pt]}

        \tabucline[2pt]{-}
    ~ &   \multicolumn{2}{c|}{\raisebox{0pt}[1\normalbaselineskip][0.5\normalbaselineskip]{$\w(v)$ if $G^D$ is of Type 1}} &  $\w(v)$ if $G^D$  is of Type 2 \\
    \tabucline[2pt]{-}
      \raisebox{0pt}[1\normalbaselineskip][0.5\normalbaselineskip]{$v \in W$} &   \multicolumn{2}{c|}{$a$}  & $a$ \\ \hline
     \multirow{2}{*}{$v \in Y_i$}  & \raisebox{0pt}[1\normalbaselineskip][0.5\normalbaselineskip]{$\textstyle a-   \frac{s-a}{i+1}$,} & if $i \ge d$ &  \multirow{2}{*}{$y_i$} \\
                  &  \raisebox{0pt}[1\normalbaselineskip][0.5\normalbaselineskip]{$a-  \frac{s-a}{d+1}$,} & if $i < d$ &  \\ \hline
      \multirow{3}{*}{$v \in B_i$}  & \raisebox{0pt}[1\normalbaselineskip][0.5\normalbaselineskip]{$a-  \frac{s-a}{i+2}- \frac{s-a}{i}$,} & if  $i > d$ & $b_{d+1}$ \qquad if $i > d$\\
                 & \raisebox{0pt}[1\normalbaselineskip][0.5\normalbaselineskip]{$a-  \frac{s-a}{d+2}-   \frac{s-a}{d+1}$,} & if  $i = d$ &   \multirow{2}{*}{$b_{i}$ \qquad \quad if $i \le d$}\\
                 & \raisebox{0pt}[1\normalbaselineskip][0.5\normalbaselineskip]{$a-2  \frac{s-a}{d+1}$,} & if  $i < d$ & \\ \hline
      \raisebox{0pt}[1\normalbaselineskip][0.5\normalbaselineskip]{$v \in R$} & \multicolumn{2}{c|}{$0$} &$0$ \\  \tabucline[2pt]{-}
\end{tabu}

\end{table}

The weight of the colored graph $G^D$ is just the sum of the weights
assigned to its vertices. Formally, $\w(G^D)=\sum_{v\in V(G)}\w(v)$.

Assume that a vertex $v\in W\cup Y$ is selected from $G^D$ in a step
of our algorithm. Hence, $v$ is recolored red in $G^{D\cup\{v\}}$.
By definition, if a neighbor $u$ of $v$ belongs to $W_i$ in $G^D$,
then $u$ is recolored yellow. Moreover, the  WY-degree of $u$
decreases by at least one, as its neighbor, $v$, was white or yellow
and now it is recolored red.
 Similarly,
if the neighbor $u$ belongs to $Y_i$ in $G^D$, then $u \in B_j$ for
a $j \le i-1$ in $G^{D\cup\{v\}}$. In the other case, if a blue
vertex $v$ is selected, $v$ is also recolored red. For any neighbor
$u$ of $v$, if $u \in W_i$ in $G^D$ then $u\in Y_j$ with $j \le i$
in $G^{D\cup\{v\}}$, and if $u \in Y_i$ in $G^D$ then $u\in B_j$
with $j \le i$ in $G^{D\cup\{v\}}$. No further vertices are
recolored, but the WY-degree of vertices from
 $N[N(v)]$ might decrease.

 Hence, assuming that the weights are nonnegative and  inequalities (1)-(8) are satisfied,
 we can observe that the weight of the colored graph and that of any vertex does not increase
 in any step of the algorithm. By conditions (1), (2), (4)-(8), the weights
 $y_i$, $b_i$, used in a  colored graph of Type 2, are not greater than the
 corresponding weights in a graph of Type~1.
 Thus, the following statement is also valid if $G^D$ belongs to
 Type~1 while $G^{D\cup\{v\}}$ belongs to Type~2.
 \begin{lem} \label{lem1}
 If the conditions (1)-(8) are satisfied, for any colored graph $G^D$ and for any vertex $v \in V(G) \setminus
 D$, the inequality $\w(G^D) \ge \w(G^{D\cup\{v\}})$ holds.
 Moreover, no vertex $u$ has greater weight in $G^{D\cup\{v\}}$ than in
 $G^D$.
 \end{lem}

\subsection{The $s$-property}
For a positive number $s$, we will say that a colored graph $G^D$
satisfies the {\it
  $s$-property}, if either $D$ is a 2-dominating set of $G$ or there exists a positive integer $k$ and a set $D^*$ of $k$
  vertices\footnote{Note that in most of the cases we will prove
  that the $s$-property holds with $|D^*|=1$. That is, we simply show
  that there exists a vertex $v$ such that the choice of $v$
  decreases $\w(G^D)$ by at least $s$.}
    such that
  $$\w(G^D)-\w(G^{D \cup D^*})\ge ks.$$

Assume that a 2-domination procedure is applied for a graph $G$
which is of order $n$. At the beginning,  we have weight $a$ on
every vertex and $\w(G^\emptyset)=an$. At the end, when $D$ is a
2-dominating set, all vertices are associated with weight 0, as they
all are contained in $R \cup B_0$. Consequently, if we show that for
every $D\subseteq V(G)$ the colored graph $G^D$ satisfies the
$s$-property, a 2-dominating set of cardinality at most $an/s$ can
be obtained, from which $\gamma_2(G) \le
  \frac{a}{s}\;n$ follows.

\begin{lem}\label{lem2}
Assume that $G$ is a graph of order $n$ and with a minimum degree of
$\delta(G)=d \ge 6$. If
 $a$, $y_0, \dots, y_{d+1}$, $b_0, \dots, b_{d+1}$ are
nonnegative numbers and $s$ is a positive number such that
conditions $(1)$--$(35)$, and for every $2\le i \le d-2$   the
inequalities $(36)$--$(41)$ are also satisfied, then for every $D
\subseteq V(G)$, the colored graph $G^D$ satisfies the $s$-property.
\end{lem}

\Proof  We  prove the lemma via a series of claims. Lemma~\ref{lem1}
will be used in nearly all argumentations here (but in most of the
cases we do not mention it explicitly). The only exception is Claim
A, which immediately follows from the definition of $s$-property.

\begin{unnumbered}{Claim A}
If $D$ is a 2-dominating set of $G$ then $G^D$ satisfies the
$s$-property.
\end{unnumbered}

\begin{unnumbered}{Claim B}
If $G^D$ belongs to Type 1, it satisfies the $s$-property.
\end{unnumbered}
\proof Let $k=\max\{i \mid W_i \cup Y_{i+1}\neq \emptyset\}$. As
$G^D$ is of Type 1, $k \ge d+1$. We assume in the next
argumentations that $G^{D \cup \{v\}}$ (or $G^{D \cup \{v'\}}$) also
is of Type 1. If this is not the case, then, by conditions (1), (2),
(4)-(8) and by the definition of the weight assignment, the decrease
in $\w(G^D)$ may be even larger than counted.

 If $W_k \neq \emptyset$, select a vertex $v\in W_k$. Each white
 neighbor $u$ of $v$ is from a class $W_i$ with $i \le k$. After the
 selection of $v$, this neighbor $u$ is recolored yellow and its
 WY-degree decreases by at least 1.\footnote{It might happen that
 the decrease is larger than 1. For example, if we have a complete
 graph $K_n$  ($n\ge 3$) with one white vertex and $n-1$ yellow
 vertices, and select the white vertex.}
   Thus, the decrease in $\w(u)$ is
 not smaller than
 $$a- \left( a-\frac{s-a}{(k-1)+1}\right)=\frac{s-a}{k}.$$
 On the other hand, each yellow neighbor $u'$ of $v$ is from a class
 $Y_{i'}$ with $i' \le k+1$. After putting $v$ into $D$, $u'$ will
 be a blue vertex with a WY-degree of at most $i'-1$. Hence, $\w(u')$ is
 decreased by at least
 $$\frac{s-a}{i'-1}\ge \frac{s-a}{k}.$$
 Since $v$ has $k$ neighbors from $W\cup Y$ in $G^D$, and the selection of $v$
 results in a decrease of $a$ in the weight of $v$, we have
 $$\w(G^D)-\w(G^{D \cup \{v\}})\ge a+ k \; \frac{s-a}{k} =s.$$
 This shows that the colored graph $G^D$ with $W_k \neq \emptyset$ satisfies the
 $s$-property.

 Now, assume that $W_k =\emptyset$. This implies
 $Y_{k+1}\neq \emptyset$ and we can select a vertex $v'\in Y_{k+1}$ in the next step of the
 procedure. As $v'$ becomes red, its weight  decreases by
 $a-\frac{s-a}{k+2}$. Each white neighbor $u$ of $v'$ has a WY-degree
 of at most $k-1$. Hence, when $u$ is recolored yellow and loses at
 least one yellow neighbor, namely $v'$, $\w(u)$ decreases by at least
 $$\frac{s-a}{(k-2)+1} >\frac{s-a}{k}.$$
 On the other hand, if $u'$ is a yellow neighbor of $v'$, we have
  the same situation as before, when a white vertex $v$ was
  put into the set $D$. That is, the decrease in $\w(u')$ is at
  least $\frac{s-a}{k}$. These imply
  $$\w(G^D)-\w(G^{D \cup \{v'\}})\ge a-\frac{s-a}{k+2}+(k+1)\frac{s-a}{k} >s$$
  and again, $G^D$ satisfies the $s$-property.
  \smallqed

  \medskip

  From now on, we consider colored graphs of Type~2. Note that the inequalities
  $$0\le y_{d+1}-b_{d} \le y_d-b_{d-1}  \le \dots \le y_2-b_1 \le
  y_1, \qquad (*)$$
  easily follow from conditions (2) and (3). Hence, if  a
  vertex $v$ is moved from $Y_i$ into $B_{i-1}$ in a step of the
  procedure, and $i \le j$ is assumed, the decrease in $\w(v)$ is at
  least $y_j-b_{j-1}$. Inequalities (1), (2) and (3) ensure similar
  estimations if $v$ is moved from $W$ into $Y_i$, from $Y_i$ into
  $B_i$, or from $B_i$ into $B_{i-1}$, and $i \le j$ is assumed.

  \begin{unnumbered}{Claim C}
If $G^D$ is a colored graph with $d-1 \le \max\{i \mid W_i \cup
Y_{i+1}\neq \emptyset\} \le d$, it satisfies the $s$-property.
\end{unnumbered}
\proof Our condition in Claim C implies that each white vertex has a
WY-degree of at most $d$ and each yellow vertex has a WY-degree of
at most $d+1$. In particular, $G^D$ is of Type~2.
  In the proof  we consider four cases.

First, assume that $W_d\neq \emptyset$ and choose a vertex $v\in
W_d$. When $v$ is put into $D$, it is recolored red and $\w(v)$
decreases by $a$. Any white neighbor $u$ of $v$ is recolored yellow
and $\deg_{WY}(u)$ decreases by at least 1. Together with condition
(1), this implies that $\w(u)$ decreases by at least $a-y_{d-1}$. A
yellow neighbor $u'$ of $v$ is recolored blue and $\deg_{WY}(u')$,
 decreases by at least 1. By $(*)$, the
weight $\w(u')$ is lowered by at least $y_{d+1}-b_d$.
  By conditions (9) and (10), $a-y_{d-1} \ge (s-a)/d$ and $y_{d+1}-b_d \ge
  (s-a)/d$. Hence,  we obtain
$$\w(G^D)-\w(G^{D \cup \{v\}})\ge a+ d \; \frac{s-a}{d}=s,$$
and $G^D$ satisfies the $s$-property.

Second, assume that $W_d= \emptyset$, but there exists a vertex
$v\in Y_{d+1}$. Let us select $v$ in the next step of the algorithm.
Then, $v$ is recolored red and $\w(v)$ decreases by $y_{d+1}$.  Each
white neighbor $u$ of $v$ has a WY-degree of at most $d-1$ in $G^D$,
and the weight  $\w(u)$ decreases by at least $a-y_{d-2}$.
Similarly,
  if  $u'$ is a yellow neighbor of  $v$,
the decrease in $\w(u')$ is not smaller than $y_{d+1}-b_d$. These
facts together with conditions (11) and (12) imply
$$\w(G^D)-\w(G^{D \cup \{v\}})\ge y_{d+1} + (d+1) \; \frac{s-y_{d+1}}{d+1}=s,$$
which proves that $G^D$ has the $s$-property.

In the third case,  $W_d \cup Y_{d+1}=\emptyset$, but there exists a
white vertex $v$ with $\deg_{WY}(v)=d-1$. Similarly to the previous
cases, but referring to  conditions (13)--(14), one can show that
$$\w(G^D)-\w(G^{D \cup \{v\}})\ge a+ (d-1) \; \frac{s-a}{d-1}=s.$$
In the last case, we assume that for each white vertex $\deg_{WY}
\le d-2$, for each yellow vertex $\deg_{WY} \le d$, and
 also that we may select a vertex $v\in Y_d$. By (15) and (16), we
 obtain
$$\w(G^D)-\w(G^{D \cup \{v\}})\ge y_{d} + d \; \frac{s-y_{d}}{d}=s.$$
This completes the proof of Claim C.   \smallqed

\medskip

\begin{unnumbered}{Claim D}
If $G^D$ is a colored graph with $\max\{i \mid W_i \cup Y_{i-1}\neq
\emptyset\} \le d-2$, and there exists a blue vertex $v$ with
$\deg_{WY}(v)\ge d+1$, then $G^D$ satisfies the $s$-property.
\end{unnumbered}
\proof Assume that $v$ is selected in the next step of the
2-domination procedure. Then,  $v$ is recolored red and $\w(v)$ is
lowered by $b_{d+1}$. Each white neighbor has a WY-degree of at most
$d-2$ and becomes yellow, while each yellow neighbor of $v$ has a
WY-degree of at most $d-1$ and becomes blue. By conditions (1) and
(2), the decrease in the weight of a white or in that of a yellow
neighbor is at least $a-y_{d-2}$ or $y_{d-1}-b_{d-1}$, respectively.
Conditions  (17) and (18) imply  $\w(G^D)-\w(G^{D \cup \{v\}})\ge
s$. \smallqed

\medskip

In the next proofs, we will use the following facts. A white vertex
does not have any red neighbors and every yellow vertex has exactly
one red neighbor. Hence, under the condition $\delta(G) \ge d$, each
white vertex $v \in W_x$ has at least $d-x$ blue neighbors, and each
$v'\in Y_y$  has at least $d-y-1$ blue neighbors. Moreover, when
this white or yellow vertex is recolored red or blue, the WY-degrees
of its $d-x$ or $d-y-1$ blue neighbors are decreased. More
precisely, if a vertex $v$ is chosen in a step of the algorithm and
$v$ is white, the sum of the WY-degrees of vertices which are blue
in $G^D$ is decreased by at least
$$d-\deg_{WY}(v)+ \sum_{w \in Y\cap N(v)} (d-1-\deg_{WY}(w)).$$
Similarly, if $v \in Y\cup B$, this decrease is at least
$$d-\deg_{WY}(v)-1+ \sum_{w \in Y\cap N(v)} (d-1-\deg_{WY}(w))$$
if $v$ is yellow, and at least
$$\deg_{WY}(v)+ \sum_{w \in Y\cap N(v)} (d-2-\deg_{WY}(w))$$
if $v$ is blue. Now, let us assume that for every blue vertex
$\deg_{WY}(u) \le j$ and for a set $B'\subseteq B$ the sum
$\sum_{u\in B'}\deg_{WY}(u)$ decreases by $z$. Then, by (2),
$\sum_{u\in B'}\w(u)$ decreases by at least $z(b_j-b_{j-1})$. This
remains valid, if for a vertex $u \in B'$, $\deg_{WY}(u)$ is reduced
by more than 1.

 \begin{unnumbered}{Claim E}
If $G^D$ is a colored graph with $d-2 \ge \max\{i \mid W_i \cup
Y_{i+1}\cup B_{i+2} \neq \emptyset\} \ge 2$, it satisfies the
$s$-property.
\end{unnumbered}
\proof Let $k=\max\{i \mid W_i \cup Y_{i+1}\cup B_{i+2} \neq
\emptyset\}$. This implies $\deg_{WY}(v) \le k$ for every white
vertex, $\deg_{WY}(v) \le k+1$ for every yellow vertex, and
$\deg_{WY}(v) \le k+2$ for every blue vertex. We consider three
cases.

If there exists a white vertex $v$ of $\deg_{WY}(v) = k$, assume
that $v$ is selected in the next step. Then, $\w(v)$ decreases by
$a$. Further, since $v$ is recolored red, the sum of the WY-degrees
of its blue neighbors decreases by at least $(d-k)$. This results in
a further change of at least $(d-k)(b_{k+2}-b_{k+1})$ in $\w(G^D)$.
If $u \in W_j$ ($j \le k$) is a white neighbor of $v$,   in
$G^{D\cup \{v\}}$ $u$ is recolored yellow and has a WY-degree of at
most $j-1$. Hence, the decrease in $\w(u)$ is at least
$$a-y_{k-1} \ge \frac{s-a -(d-k)(b_{k+2}-b_{k+1})}{k},$$
where the last inequality follows from (36) substituting $i=k$.
Consider now a yellow neighbor $u'$ of $v$. After the choice of $v$,
$u'$ is recolored blue and $\w(u')$ decreases by at least
$y_{k+1}-b_k$. Taking into account the decreases in the weights of
vertices from $N(u')\cap B$, the recoloring of each such $u'$
contributes to the  decrease of $\w(G^D)$ with at least
$$y_{k+1}-b_k+(d-(k+1)-1)(b_{k+2}-b_{k+1}) \ge  \frac{s-a -(d-k)(b_{k+2}-b_{k+1})}{k},$$
where the lower bound follows from (37) substituting $i=k$.
Therefore, if $v \in W_k$,
$$\w(G^D)-\w(G^{D \cup \{v\}})\ge a+ (d-k)(b_{k+2}-b_{k+1})+ k \; \frac{s-a -(d-k)(b_{k+2}-b_{k+1})}{k}=s.$$
Consequently,  $G^D$ has the $s$-property if $W_k \neq \emptyset$.

In the following two cases, we  count $\w(G^D)-\w(G^{D \cup \{v\}})$
in a similar way. Assume that $W_k = \emptyset$ but $Y_{k+1} \neq
\emptyset$, and choose a vertex $v$ from $Y_{k+1}$.
 Vertex $v$ is recolored red and the WY-degrees of its blue neighbors
 decrease. This contributes to the difference $\w(G^D)-\w(G^{D \cup \{v\}})$ with at least $y_{k+1} +
 (d-k-2)(b_{k+2}-b_{k+1})$.
 Further, if $u$ is a white neighbor of $v$ then $\deg_{WY}(u) \le k-1$ in $G^D$. Once $v$ is recolored
 red, $\w(u)$ decreases by at least $a-y_{k-2}$. By condition (38), it is
 not smaller than
 $(s-y_{k+1}-(d-k-2)(b_{k+2}-b_{k+1}))/(k+1)$. If $u'$ is a yellow
 neighbor of $v$, then $u'$ will be blue in $G^{D \cup \{v\}}$ and  the
 WY-degrees in $B\cap N(u')$ are decreased. Consequently, and also referring to (39),
  each yellow neighbor $u'$ contributes to the decrease of $\w(G^D)$ with at least
 $$y_{k+1}-b_k +(d-k-2)(b_{k+2}-b_{k+1})\ge
 \frac{s-y_{k+1}-(d-k-2)(b_{k+2}-b_{k+1})}{k+1}.$$
 In total, $v$ has $k+1$ neighbors from $W\cup Y$, and we have
 $$\w(G^D)-\w(G^{D \cup \{v\}})\ge y_{k+1} +
 (d-k-2)(b_{k+2}-b_{k+1})+(k+1) \; \frac{s-y_{k+1}-(d-k-2)(b_{k+2}-b_{k+1})}{k+1}=s,$$
which proves that $G^D$ satisfies the $s$-property.

In the third case, $W_k \cup Y_{k+1} =\emptyset$ and we have a blue
vertex $v$ with $\deg_{WY}(v)=k+2$. Selecting $v$ in the next step
of the procedure, $v$ will be recolored red and $\w(v)$ becomes 0.
Each white neighbor $u$ of $v$ is recolored yellow and has a
decrease of at least $a-y_{k-1}$ in $\w(u)$ (in this case,
$\deg_{WY}(u)$ might be unchanged). Moreover, each yellow neighbor
$u'$ of $v$ is recolored blue and the weights of the vertices from
$N(u')\cap B$ are also decreased. Then, the recoloring of $u'$
contributes to the decrease of $\w(G^D)$ by at least
$$y_k-b_k +(d-k-2)(b_{k+2}-b_{k+1}) \ge \frac{s-b_{k+2}}{k+2},$$
where the inequality follows from (40). On the other hand, by (41),
we have $a-y_{k-1} \ge (s-b_{k+2})/(k+2)$. We may conclude that
$$\w(G^D)-\w(G^{D \cup \{v\}})\ge b_{k+2}+ (k+2) \;
\frac{s-b_{k+2}}{k+2} = s.$$
 Thus, in the third case $G^D$ also satisfies the $s$-property.
 \smallqed

 \begin{unnumbered}{Claim F}
Let $G^D$ be a colored graph with $  \max\{i \mid W_i \cup
Y_{i+1}\cup B_{i+2} \neq \emptyset\} = 1$ such that there exists an
edge between $W$ and $Y$. Then, $G^D$ satisfies the $s$-property.
\end{unnumbered}
\proof Choose a white vertex $v$ whose only neighbor from $W\cup Y$
is a yellow vertex $u$  in $G^D$. By our condition, $\deg_{WY}(u)
\le 2$. In $G^{D \cup \{v\}}$, the vertex $v$ is recolored red and
$u \in B_1$. Moreover, in $G^D$, $v$ and $u$ has at least $d-1$ and
$d-3$ blue neighbors, respectively. By condition (19),
$$\w(G^D)-\w(G^{D \cup \{v\}})\ge
a+(d-1)(b_3-b_2)+(y_2-b_1)+(d-3)(b_3-b_2)\ge s,$$ and $G^D$ has the
$s$-property. \smallqed

Henceforth, we may assume that there are no edges between $W$ and
$Y$.

\begin{unnumbered}{Claim G}
If $G^D$ is a colored graph with $  \max\{i \mid W_i \cup
Y_{i+1}\cup B_{i+2} \neq \emptyset\} = 1$ and $Y_2\neq \emptyset$,
then $G^D$ has the $s$-property.
\end{unnumbered}

\proof Consider a vertex $v \in Y_2$ in $G^D$. As supposed,
 it has no white neighbors. Hence, $v$ is adjacent to two
vertices, say $u_1$ and $u_2$, which are from $Y_2 \cup Y_1$. Then,
in $G^{D \cup \{v\}}$, $v$ is recolored red, $u_1$ and $u_2$ are
recolored blue and belong to $B_1 \cup B_0$. The decrease in
$\sum_{w \in B\cap (N(v)\cup N(u_1) \cup N(v_2))}\deg_{WY}(w)$ is at
least $3(d-3)$. Then, also using (20),
$$\w(G^D)-\w(G^{D \cup \{v\}})\ge
y_2+2(y_2-b_1)+3(d-3)(b_3-b_2) \ge s.$$
 This proves the claim. \smallqed

 \begin{unnumbered}{Claim H}
 If $G^D$ is a colored graph with $
\max\{i \mid W_i \cup Y_{i}\cup B_{i+2} \neq \emptyset\} = 1$, it
satisfies the $s$-property.
\end{unnumbered}
\proof Suppose for a contradiction that there exits a colored graph
$G^D$ which satisfies the condition of our claim but does not have
the $s$-property. First, let us assume   $W_1\neq \emptyset$ and
recall that each white vertex with $\deg_{WY}(v)=1$ has a white
neighbor of the same type.  We consider the following cases:
\begin{itemize}
\item[$(i)$] If there exists a vertex $v_1 \in W_1$ with a white neighbor
$v_2$, and with a blue neighbor $u$ from $B_3$ such that $u$ is not
adjacent to $v_2$, we assume that in two consecutive steps of the
procedure  $v_2$ and $u$ are chosen. Then, $v_2$ and $u$ are
recolored red, and $v_1$ becomes blue with a WY-degree of 0. This
contributes to the decrease of $\w(G^D)$ with $2a+b_3$. The total
weight of the further blue neighbors of $v_1$ and $v_2$ decreases by
at least $((d-2)+(d-1))(b_3-b_2)$. If $u$ has a white neighbor $w$
in $G^D$, $w$ becomes yellow and contributes to the decrease of
$\w(G^D)$ with at least $a-y_1$. By (21), it  is  not smaller than
$(2s-2a-b_3-(2d-3)(b_3-b_2))/2$. If $w'$ is a yellow neighbor of $u$
in $G^D$, then it is recolored blue and $\deg{WY}(w')$ is either 1
or 0 in $G^{D\cup \{v_2,u\}}$. Further, the weights of the at least
$d-3$ blue neighbors of $w'$ which are different from $u$ are also
decreased. In total, $w'$ contributes to the decrease of $\w(G^D)$
with at least
$$y_1-b_1+(d-3)(b_3-b_2) \ge \frac{2s-2a-b_3-(2d-3)(b_3-b_2)}{2},$$
where the last inequality is equivalent to (22). Therefore, we have
$$\w(G^D)-\w(G^{D \cup \{v_2,u\}})\ge
2a+b_3+(2d-3)(b_3-b_2)+ 2 \;
\frac{2s-2a-b_3-(2d-3)(b_3-b_2)}{2}=2s,$$
 and the $s$-property  would be satisfied by $G^D$. This contradicts
 our assumption.

 \item[$(ii)$] Since $G^D$ is supposed to be a counterexample,
  if a  blue vertex $u \in
 B_3$ is adjacent to a white vertex then it is also adjacent to the
 white neighbor of it. If we have two adjacent white vertices $v_1$
 and $v_2$ which have only one (common) neighbor $u$ from $B_3$, choose
 $v_1$ and $u$ in the next two steps of the procedure. Then, $v_1$
 and $u$ are recolored red, while $v_2$ is recolored blue and has a
 WY-degree of 0. Their weights are decreased by $2a+b_3$. All the
 further blue neighbors of $v_1$ and $v_2$ belong to $B_2 \cup B_1$
 in $G^D$. The WY-degrees of these blue vertices are reduced, which
 contributes to the difference $\w(G^D)-\w(G^{D \cup \{v_1,u\}})$
 with  at least $2(d-2)(b_2-b_1)$. The blue vertex $u$ has one white
 or yellow neighbor $w$ which is different from $v_1$ and $v_2$. If
 $w$ is white,  it is from $W_0$, as otherwise $w$, its white neighbor, and
 $u$ would satisfy the assumption in case $(i)$. Hence, when $w$ is recolored yellow, $\w(w)$
 decreases by $a-y_0$, and
 $$\w(G^D)-\w(G^{D \cup \{v_1,u\}})\ge
2a+b_3+2(d-2)(b_2-b_1)+ a-y_0, $$ which is at least $2s$ by
condition (23).
  If $w$ is yellow then $w \in Y_1 \cup Y_0$. When $w$ is recolored
  blue, the WY-degrees of its blue neighbors are also reduced.
  These contribute to the difference $\w(G^D)-\w(G^{D \cup
  \{v_1,u\}})$ with at least $y_1-b_1+ (d-3)(b_3-b_2)$.
  Therefore, referring to (24),
  $$\w(G^D)-\w(G^{D \cup \{v_1,u\}})\ge
2a+b_3+2(d-2)(b_2-b_1)+ y_1-b_1+ (d-3)(b_3-b_2)\ge 2s. $$
  We infer that in the counterexample $G^D$ we cannot have  a white vertex
  in $W_1$ that has exactly one neighbor from
 $B_3$.

  \item[$(iii)$] Now assume that $v_1, v_2 \in W_1$ and their
  neighbors $u_1$ and $u_2$ are from $B_3$ in $G^D$. Choose $u_1$
  and  $u_2$ and consider $G^{D \cup \{u_1,u_2\}}$. Here, $v_1$ and
  $v_2$ are blue vertices of WY-degree 0, while $u_1$ and $u_2$ are
  red. In $G^D$, each blue neighbor of $v_1$ and $v_2$ which is
  different from $u_1$ and $u_2$ is either from $B_2$ or it is a
  further common neighbor of $v_1$ and $v_2$ from $B_3$. In the worst case, the
  decrease in their weights contributes to $\w(G^D)-\w(G^{D \cup
  \{u_1,u_2\}})$ with $(d-3)(b_3-b_1)$. Finally, $u_1$ and $u_2$
  have neighbors from $W_0 \cup Y_1 \cup Y_0$. It is enough to
  consider the following cases.
  \begin{itemize}
  \item $u_1$ and $u_2$ have a common neighbor $w \in W_0$. Then,
  $w$ is recolored blue.
 The  weight of $w$ and that of its blue neighbors (different from $u_1$ and $u_2$) decrease
  by at least
  $a+(d-2)(b_3-b_2)$.
  Then, by (25) and by our earlier observations
  $$\w(G^D)-\w(G^{D \cup \{u_1,u_2\}}) \ge
   2a+2b_3+(d-3)(b_3-b_1)+a+(d-2)(b_3-b_2) \ge 2s.$$
   Hence, in a counterexample we cannot have this case.
   \item $u_1$ and $u_2$ have a common neighbor $w \in Y_1$. Then,
  $w$ is recolored blue and moved to $B_1$ in $G^{D \cup
  \{u_1,u_2\}}$. Also, the  weights of its blue neighbors
  decrease. These contribute to the difference $\w(G^D)-\w(G^{D \cup
  \{u_1,u_2\}})$ with at least
  $y_1-b_1+ (d-4)(b_3-b_2)$, and we have
   $$\w(G^D)-\w(G^{D \cup \{u_1,u_2\}}) \ge
   2a+2b_3+(d-3)(b_3-b_1)+y_1-b_1+ (d-4)(b_3-b_2) \ge 2s,$$ where the last inequality follows from
   (26). Again, this case is not possible in a counterexample.
   \item $u_1$ and $u_2$ have two different neighbors, namely  $w_1$ and $w_2$, from
   $W_0$. Then, $w_1$ and $w_2$ are recolored yellow and we have
  \begin{align*}\w(G^D)-\w(G^{D \cup \{u_1,u_2\}}) &\ge
   2a+2b_3+(d-3)(b_3-b_1)+2(a-y_0)\\
    &\ge 2a +b_3 +3(b_3-b_2)+2(d-3)(b_3-b_2)+2(a-y_1) \ge 2s.
   \end{align*}
  Here, we used (21) and the inequalities $b_3\ge 3(b_3-b_2)$ and $b_3-b_1
  \ge 2(b_3-b_2)$ which follow from (2).
  \end{itemize}
  \end{itemize}
 We have shown that there are no edges between $W_1$ and $B_3$ if
 $G^D$ is a counterexample to Claim F. In what follows we prove that
 $B_3=\emptyset$ and $Y_1=\emptyset$.

 Suppose that $B_3\neq  \emptyset$ and  choose a vertex $v$
 from $B_3$. As it has been shown, all white and yellow neighbors of $v$ belong to $W_0 \cup Y_1 \cup
 Y_0$. If $u$ is a white neighbor, $\w(u)$ decreases by $a-y_0$, and
 if $u'$ is yellow, its recoloring contributes to the decrease of $G^D$ by at
 least $y_1-b_1+ (d-3)(b_3-b_2)$. By conditions (27) and (28),
 $$\w(G^D)-\w(G^{D \cup \{v\}}) \ge b_3+ 3 \; \frac{s-b_3}{3}=s.$$
 Hence, in the counterexample each blue vertex is of a WY-degree of at
 most 2.

 Suppose now that $Y_1 \neq \emptyset$ and choose a vertex $v$ from it.
 Since $v$ cannot have a neighbor from $W$, it must have a neighbor
 $u$ from $Y_1$. In $G^{D \cup \{v\}}$, $v$ is recolored red, $u$ is
 recolored blue with a WY-degree 0, and each of their at least
 $2(d-2)$ blue neighbors has a decrease of at least $b_2-b_1$ in its weight. Hence, we have
 $$\w(G^D)-\w(G^{D \cup \{v\}}) \ge 2y_1+2(d-2)(b_2-b_1),$$
 which is at least $s$ by (29). We may conclude that $Y_1=\emptyset$
 holds in our counterexample.

 Assume that $W_1$ is not empty. Then, $W_1$ consists of pairs of adjacent
 vertices, we refer to which as ``white pairs".

 First, suppose that there exits a white pair $v_1,v_2$ and a vertex $u\in B_2$
 such that $u$ is adjacent to $v_1$ and nonadjacent to $v_2$. In the
 next two steps of the procedure we choose $v_2$ and $u$. Then,
 $v_2$ and $u$ are recolored red, $v_1$ becomes a blue vertex of
 WY-degree 0. The WY-degrees of blue neighbors of $v_1$ and $v_2$
 are also reduced. In total, these result in a decrease of at least $2a+b_2+(2d-3)(b_2-b_1)$
 in $\w(G^D)$. Moreover, $u$ has a white or a yellow neighbor $w$
 different from $v_1$. For the cases  $w\in  W_1$ and $w\in Y_1\cup
 Y_0$ we have the following inequalities by (30) and (31),
 respectively.
 $$\w(G^D)-\w(G^{D \cup \{v_2,u\}}) \ge 2a+b_2+ (2d-3)(b_2-b_1)+
 (a-y_1) \ge 2s$$
 $$\w(G^D)-\w(G^{D \cup \{v_2,u\}}) \ge 2a+b_2+ (2d-3)(b_2-b_1)+
 (y_1-b_1+(d-3)(b_2-b_1)) \ge 2s$$
 We may infer that $G^D$ has the $s$-property, which is a contradiction.
 Hence, if a blue vertex from $B_2$ is adjacent to a vertex from $W_1$, then it is also adjacent to the
 other vertex from that white pair.

Now, consider any white pair $v_1,v_2$ and choose these two vertices
in two consecutive steps of the procedure. As a result, $v_1$ and
$v_2$ are recolored red and all their blue neighbors are of
WY-degree 0. Since $b_2-b_1 \le b_1-b_0=b_1$, the worst case is when
$v_1$ and $v_2$ share $d-1$ blue neighbors from $B_2$ in $G^D$. By
(32), we have
 $$\w(G^D)-\w(G^{D \cup \{v_1,v_2\}}) \ge 2a+(d-1)b_2 \ge 2s,$$
 contradicting our assumption that $G^D$ is  counterexample.

 Consequently, if $\max\{i \mid W_i \cup Y_{i}\cup B_{i+2} \neq \emptyset\} =
 1$ then $G^D$  has the $s$-property, as stated in Claim H.
 \smallqed

 What remains to consider after Claims A-H is  the case when $D$ is not a
 2-dominating set that is $W\cup Y \neq \emptyset$ but all white and
 yellow vertices are of WY-degree 0 and all blue vertices have
 a WY-degree of at most 2.

 First, suppose that we have an edge between $B_2$ and $Y_0$. Then,
 choose a blue vertex $v\in B_2$ which has a yellow neighbor $u$.
Vertex $v$ has a further neighbor $u'$ from $W_0 \cup Y_0$.
Depending on the color of $u'$, we can use either (33) or (34) and
obtain the following inequalities. If $u'$ is yellow,
$$\w(G^D)-\w(G^{D \cup \{v\}}) \ge b_2+2y_0+2(d-2)(b_2-b_1) \ge s.$$
If $u'$ is white
$$\w(G^D)-\w(G^{D \cup \{v\}}) \ge b_2+y_0+(d-2)(b_2-b_1) +a-y_0 \ge s.$$
Thus, in these cases  $G^D$ has the $s$-property.

Now assume that $Y_0\neq \emptyset$ and choose a vertex $v$ from
$Y_0$. We have just shown that $v$ has no neighbors from $B_2$.
Hence, $v$ has at least $d-1$ blue neighbors from $B_1$. Together
with (35), these imply
$$\w(G^D)-\w(G^{D \cup \{v\}}) \ge y_0+(d-1)b_1\ge s,$$
and $G^D$ has the $s$-property.

Finally, we assume that $Y=\emptyset$, but we have $x$ vertices in
$W_0$, $z_2$  vertices in $B_2$ and $z_1$ vertices in $B_1$, Thus,
$\w(G^D)=xa + z_2b_2+ z_1b_1$. On the other hand, counting the
number of edges between $W_0$ and $B_2\cup B_1$ in two different
ways, $dx \le 2z_2 +z_1.$ Consider $G^{D\cup Y_0}$, that is assume
that in $x$ consecutive steps we select all white vertices. Clearly,
in $G^{D\cup Y_0}$ every vertex has a weight of 0. Hence,
\begin{align*}
 \w(G^D)-\w(G^{D \cup Y_0}) &=xa + z_2b_2+ z_1b_1 \ge xa+
 (2z_2+z_1)\min \left\{\frac{b_2}{2}, b_1\right\}\\
 &\ge  xa+dx \; \frac{b_2}{2} \ge xs.
\end{align*}
The last inequality is a consequence of (32), and $b_2/2 \le b_1$
follows from $b_2-b_1 \le b_1$.

The cases discussed in our proof together cover all possibilities,
hence every colored graph $G^D$ satisfies the $s$-property under the
conditions of Lemma~\ref{lem2}. \qed

\medskip

As we discussed it at the beginning of this section,
Theorem~\ref{th1} is an immediate consequence of Lemma~\ref{lem2}.

\section{Concluding remarks}

Finally, we make some remarks on the algorithmic aspects of our
proof. In Table~\ref{Atable11}, we compared the upper bounds
obtained  by our Theorem~\ref{th1} and those proved in \cite{HV}
with probabilistic method. Our upper bounds on $\gamma_2(G)$ improve
the earlier best results if the minimum degree $\delta$ is between 6
and 21. Nevertheless the
 algorithm, which is behind  our proof, can also be useful for
 $\delta \ge 22$, as we can guarantee the determination of a
 2-dominating set of bounded size for each input graph.

We can identify two  different algorithms based on the proof in
Section~\ref{sec-proof}. For the first version, we do not need to
count the weights assigned to the vertices. We just consider the
 list of instructions below  and in each step of the algorithm we follow
the first one which is applicable.
\begin{enumerate}
\item If $k=\max\{i \mid W_i \cup Y_{i+1}\neq \emptyset\} \ge d-1$
and $W_k\neq \emptyset$, choose a vertex from $W_k$.
\item If $k=\max\{i \mid W_i \cup Y_{i+1}\neq \emptyset\} \ge d-1$,
 choose a vertex from $Y_{k+1}$.
 \item  If $k=\max\{i \mid B_{i}\neq \emptyset\} \ge d+1,$
  choose a vertex from $B_k$.
  \item If $2\le k=\max\{i \mid W_i \cup Y_{i+1}\cup B_{i+2}\neq \emptyset\} \le d-2$
and $W_k\neq \emptyset$, choose a vertex from $W_k$.
\item If $2\le k=\max\{i \mid W_i \cup Y_{i+1}\cup B_{i+2}\neq \emptyset\} \le d-2$
and $Y_{k+1}\neq \emptyset$, choose a vertex from $Y_{k+1}$.
\item If $2\le k=\max\{i \mid W_i \cup Y_{i+1}\cup B_{i+2}\neq \emptyset\} \le d-2,$
 choose a vertex from $B_{k+2}$.
 \item If there exists a white vertex $v$ with a yellow neighbor,
 choose $v$.
 \item If $Y_2 \neq \emptyset$, choose a vertex from it.
 \item If there exist two adjacent white vertices $v_1$ and $v_2$
 such that $v_1$ has a neighbor $u$ from $B_3$ which is not adjacent
 to $v_2$, choose $v_2$ and $u$.
 \item If there exists a vertex $v$ in $W_1$, which has exactly one
 neighbor, say $u$, in $B_3$, choose $v$ and $u$.
 \item If there exists a vertex $v$ in $W_1$, which has at least two
 neighbors in $B_3$, choose two vertices from $N(v) \cap B_3$.
 \item If $B_3\neq \emptyset$, choose a vertex from it.
 \item If $Y_1\neq \emptyset$, choose a vertex from it.
 \item If there exist two adjacent white vertices $v_1$ and $v_2$
 such that $v_1$ has a neighbor $u$ from $B_2$ which is not adjacent
 to $v_2$, choose $v_2$ and $u$.
 \item If there exist two adjacent white vertices, choose such two
 vertices.
 \item If there exists a blue vertex $v \in B_2$ which has at least
 one yellow neighbor, choose $v$.
 \item If $Y\neq \emptyset$, choose a yellow vertex.
 \item Choose all the white vertices.
\end{enumerate}

By a slightly different interpretation, we can define a 2-domination
algorithm based on the weight assignment introduced in
Section~\ref{sec-proof}. Then, in each step, we choose a vertex $v$
such that the decrease $\w(G^D)-\w(G^{D\cup \{v\}})$ is the possible
largest. The exceptions are those steps where $G^D$ would be treated
by instructions
 9, 10, 11, 14, 15 or 18 of the previous algorithm. In these cases,
 the greedy choice concerns the maximum decrease of $\w(G^D)$ in
 two (or more) consecutive steps.

\medskip
\section*{Acknowledgements}
Research of Csilla Bujt\' as  was supported by the National
Research, Development and Innovation Office -- NKFIH under the grant
SNN 116095.

\end{document}